\documentclass[12pt]{amsart}

\newtheorem{thm}{Theorem}
\newtheorem{prop}[thm]{Proposition}
\newtheorem{lem}[thm]{Lemma}
\newtheorem{cor}[thm]{Corollary}

\theoremstyle{remark}
\newtheorem{rem}[thm]{Remark}

\theoremstyle{definition}
\newtheorem{defn}[thm]{Definition}
\newtheorem{ex}[thm]{Example}

\newcommand{\C}{\mathbb{ C}}

\newcommand{\R}{\mathbb{ R}}
\newcommand{\HH}{\mathbb{ H}}

\title{On products of harmonic forms}
\author{D.~Kotschick}
\address{Mathematisches Institut, Universit\"at M\"unchen,
Theresienstr.~39, 80333 M\"unchen, Germany}
\email{dieter@member.ams.org}
\date{\today; MSC2000: 53C25, 57R57, 58A14, 57R17}

\begin{document}

\begin{abstract}
We prove that manifolds admitting a Riemannian metric for
which products of harmonic forms are harmonic satisfy strong
topological restrictions, some of which are akin to properties 
of flat manifolds. Others are more subtle, and are related to
symplectic geometry and Seiberg-Witten theory.

We also prove that a manifold admits a metric with harmonic
forms whose product is not harmonic if and only if it is not a 
rational homology sphere.
\end{abstract}

\maketitle



\section{Introduction}\label{s:intro}

On a general Riemannian manifold, wedge products of harmonic
forms are not usually harmonic. But there are some examples
where this does happen, like compact globally symmetric spaces.
For these the harmonic forms coincide with the invariant ones,
and the latter are clearly closed under products. See~\cite{DFN},
pp.~10-13. 

Sullivan~\cite{tokyo} observed that ``There are topological 
obstructions for $M$ to admit a metric in which the product of
harmonic forms is harmonic.'' The reason Sullivan gave is that
if the product of harmonic forms is harmonic, then the rational
homotopy type is a formal consequence of the cohomology ring.
Therefore, manifolds which are not formal in this sense cannot
admit a metric for which the products of harmonic forms are
harmonic.

This motivates the following:
\begin{defn}
A Riemannian metric is called (metrically) formal if all wedge 
products of harmonic forms are harmonic.

A closed manifold is called geometrically formal if it admits
a formal Riemannian metric. 
\end{defn}

Thus geometric formality implies formality in the sense of
Sullivan. Compact globally symmetric spaces are metrically formal, 
as are arbitrary Riemannian metrics on rational homology spheres.
Further examples can be generated by taking products, because
the product of two formal metrics is again formal.

In this paper we describe a number of elementary topological 
obstructions for geometric formality of closed oriented manifolds.
These obstructions are independent of formality in the sense of 
rational homotopy theory, and are often nonzero on formal manifolds. 
The simplest obstruction is the product of the first Betti number 
and the Euler characteristic. In small dimensions these elementary
obstructions are strong enough to imply:

\begin{thm}\label{t:main}
If $M$ is a closed oriented geometrically formal manifold of 
dimension $\leq 4$, then $M$ has the real cohomology algebra
of a compact globally symmetric space.
\end{thm}
It is however not true that $M$ is a globally symmetric space,
even up to homotopy. We give many examples of this in dimensions
$3$ and $4$. 

We also give examples of $4$-manifolds which do have the real 
cohomology algebra of a compact symmetric space, but are not 
geometrically formal. This is detected by some more subtle 
obstructions coming from symplectic geometry and Seiberg-Witten 
gauge theory. 

The pattern of the arguments presented here is that metric
formality is a weakening of a reduction of holonomy. For 
example, it implies that harmonic forms have constant length,
though it does not imply that they are parallel. Nevertheless,
the more harmonic forms there are, the stronger the constraints.

In Section~\ref{s:generic} we prove that every manifold which 
is not a rational homology sphere admits a non-formal Riemannian 
metric.

This paper was motivated by an example pointed out by D.~Toledo,
which we describe in Section~\ref{s:two}, and which is 
related to joint work in progress of the author with H.~Endo.
Further impetus came from a question posed by
D.~Huybrechts and U.~Semmelmann concerning products of 
harmonic forms on Calabi-Yau manifolds, see~\cite{Huy}.

\section{Motivation: dimension two}\label{s:two}

Let us consider first the case of a closed oriented surface $\Sigma$.
If its genus is $0$ or $1$, then there are globally symmetric Riemannian
metrics. If the genus of $\Sigma$ is $\geq 2$, there are nontrivial 
harmonic $1$-forms for all metrics, but every $1$-form has zeros. 
In this case every wedge product of $1$-forms also has zeros, but for 
cohomological reasons cannot vanish identically in all cases. The only 
harmonic $2$-forms are the constant multiples of the Riemannian volume 
form, so there cannot be any formal Riemannian metric. This proves
Theorem~\ref{t:main} in the $2$-dimensional case.

The argument above, pointed out to the author by D.~Toledo, shows that 
harmonicity of products of harmonic forms can fail on compact locally
rather than globally symmetric spaces, contradicting a statement 
in~\cite{GM}, p.~158. 

Note that $\Sigma$ is formal in the sense of Sullivan, but that the 
nonvanishing of $b_1(\Sigma ) \cdot \chi (\Sigma )$ obstructs 
geometric formality.   

On the sphere every metric is formal, because there are no interesting 
harmonic forms. However, when there are enough harmonic forms,
the harmonicity of their products is a restriction on the metric
enforcing rigidity.
 
\begin{thm}\label{t:two}
Every formal Riemannian metric on the two-torus is flat.
\end{thm}
\begin{proof}
Let $g$ be a formal metric on $T^2$, and $\alpha$ a nontrivial harmonic 
$1$-form. Then $*\alpha$ and
\begin{equation}\label{eq:first}
\alpha\wedge *\alpha = \vert\alpha\vert^2 dvol_g
\end{equation}
are also harmonic and so $\alpha$ has constant length. In particular
it has no zeros. It follows that
$\alpha (p)$ and $*\alpha (p)$ span $T_p\Sigma$ for all $p\in\Sigma$.
As $\vert a\vert$ is constant for every constant linear combination
$a$ of $\alpha$ and $*\alpha$, the Bochner formula
\begin{equation}\label{eq:Bochner}
\Delta (a) = \nabla^*\nabla a +Ric (a)
\end{equation}
allows us to compute:
$$
0=\Delta (\frac{1}{2}\vert a\vert^2)= g(\nabla^*\nabla a,a)-
\vert\nabla a\vert^2 =-g(Ric (a),a)-\vert\nabla a\vert^2 \ .
$$
This shows that the (Ricci) curvature is everywhere nonpositive. But by
the Gau{\ss}-Bonnet theorem this implies that $g$ is flat. 
\end{proof}

\section{Elementary obstructions}\label{s:elem}

Let $M$ be a closed oriented manifold of dimension $n$, and $g$
a formal Riemannian metric on $M$. As usual, we extend $g$ to 
spaces of differential forms.

\begin{lem}\label{l:length}
The inner product of any two harmonic $k$-forms is a constant
function. In particular, the length of any harmonic form is constant.
\end{lem}
\begin{proof}
That the length of any harmonic form is constant follows from 
equation~\eqref{eq:first}. The more general statement follows
by polarisation.
\end{proof}

\begin{lem}\label{l:comb}
Suppose $\alpha_1,\ldots,\alpha_m$ are orthogonal harmonic $k$-forms.
Then 
$$
\alpha = \sum_{i=1}^m f_i\alpha_i
$$
is harmonic if and only if the functions $f_i$ are all constant.
\end{lem}
\begin{proof}
If $\alpha$ is harmonic, then $g(\alpha,\alpha_i)=f_i \vert\alpha_i\vert^2$ 
is constant by Lemma~\ref{l:length}. Using that the length of 
$\alpha_i$ is also constant by Lemma~\ref{l:length}, we conclude that
$f_i$ is constant.

The converse is trivial. 
\end{proof}

Lemma~\ref{l:length} implies that harmonic forms which are linearly
independent at some point are linearly independent everywhere.
Systems of linearly independent harmonic forms can be orthonormalised
using constant coefficients.

We can now generalise the discussion in Section~\ref{s:two} to higher dimensions.
\begin{thm}\label{t:elem}
Suppose the closed oriented manifold $M^n$ is geometrically formal. 
Then
\begin{enumerate}
\item the real Betti numbers of $M$ are bounded by $b_k(M)\leq b_k(T^n)$, and
\item if $n=4m$, then $b_{2m}^{\pm}(M)\leq b_{2m}^{\pm}(T^n)$.
\item The first Betti number $b_1(M)\neq n-1$.
\end{enumerate}
\end{thm}
\begin{proof}
Fix a formal Riemannian metric on $M$.
It follows from the above Lemmas that the number of linearly
independent harmonic $k$-forms is at most the rank of the vector
bundle $\Lambda^k$. Similarly, when the dimension is $4m$,
the number of self-dual or anti-self-dual harmonic forms in 
the middle dimension is bounded by the rank of $\Lambda^{2m}_{\pm}$.

Suppose now that $\alpha_1,\ldots,\alpha_{n-1}$ are linearly
independent harmonic $1$-forms. Then 
$*(\alpha_1\wedge\ldots\wedge\alpha_{n-1})$ is also a harmonic
$1$-form, and is linearly independent of $\alpha_1,\ldots,\alpha_{n-1}$.
Thus $b_1(M)\geq n-1$ implies $b_1(M)= n$.
\end{proof}

There is an uncanny similarity here with the classification of 
flat Riemannian manifolds~\cite{wagner}, which satisfy all the 
conclusions of Theorem~\ref{t:elem}. We can push this further:
\begin{thm}\label{t:fiber}
Suppose the closed oriented manifold $M^n$ is geometrically
formal. 
If $b_1(M)=k$, then there is a smooth submersion $\pi\colon M
\rightarrow T^k$, for which $\pi^*$ is an injection of  
cohomology algebras. In particular, if $b_1(M)=n$, then $M$
is diffeomorphic to $T^n$. 
In this case every formal Riemannian metric is flat.
\end{thm}
\begin{proof}
Fix a formal Riemannian metric $g$ on $M$.
We consider the Albanese or Jacobi map $\pi\colon M\rightarrow T^k$ 
given by integration of harmonic $1$-forms. As the harmonic $1$-forms
have constant lengths and inner products, $\pi$ is a submersion.
It induces an isomorphism on $H^1$, and products of linearly
independent harmonic $1$-forms are never zero, but are harmonic
because the metric is formal.

In the case $b_1(M)=n$, we conclude that $M$ is a covering of 
$T^n$, and is therefore a torus itself. Every formal metric
on $T^n$ must be flat because it admits an orthonormal framing 
by harmonic $1$-forms.
\end{proof}

In the case $b_1(M)=1$ we have a partial converse to 
Theorem~\ref{t:fiber}:
\begin{thm}\label{t:fiberconv}
Let $M$ be a closed oriented $n$-manifold which fibers smoothly
over $S^1$. If $b_1(M)=1$ and $b_k(M)=0$ for $1<k<n-1$, then
$M$ is geometrically formal.
\end{thm}
\begin{proof}
Suppose $M$ fibers over $S^1$ with fiber $F$ and monodromy 
diffeomorphism $\phi\colon F\rightarrow F$. By Moser's Lemma, 
we may assume that $\phi$ preserves a volume form $\epsilon$ on 
$F$, so that its pullback to $F\times\R$ descends to $M$ as a 
closed form which is a volume form along the fibers. We can find 
a Riemannian metric on $M$ for which $*\epsilon =\alpha$ is the 
closed $1$-form defining the fibration over $S^1$, and has constant 
length. Then $\alpha$ and $\epsilon$ generate the harmonic forms 
in degree $1$ and $n-1$, and their product is harmonic.
\end{proof}

Flat manifolds satisfy further topological restrictions, for
example their Euler characteristics vanish. In our present
context we have:
\begin{thm}\label{t:Euler}
Suppose the closed oriented manifold $M^n$ is geometrically formal. 
\begin{enumerate}
\item If $b_k(M)\neq 0$, then $e(\Lambda^k)=0$, and
\item if $n=4m$ and $b_{2m}^{\pm}(M)\neq 0$, then 
$e(\Lambda^{2m}_{\pm})=0$.
\end{enumerate}
In particular the Euler characteristic of $M$ vanishes if 
$b_1(M)\neq 0$.
\end{thm}
\begin{proof}
This follows from the obstruction-theory definition of the 
Euler class, and the fact that every nontrivial 
harmonic form has no zeros because of~\eqref{eq:first}.
\end{proof}

\section{Dimension three}\label{s:three}

If $M$ is a closed oriented geometrically formal $3$-manifold, 
then by Theorem~\ref{t:elem} we have $b_1(M)\in\{ 0,1,3\}$.
If the first Betti number is maximal, then Theorem~\ref{t:fiber} 
says that $M$ is the 3-torus.
At the other extreme, if the first Betti number is zero, then 
$M$ is a rational homology sphere. Clearly every metric on every
such manifold is formal.

Thus, the only interesting case is that of first Betti number one.
Then the real cohomology algebra is that of the globally symmetric
space $S^2\times S^1$, so that Theorem~\ref{t:main} is proved in the 
$3$-dimensional case.

Theorems~\ref{t:fiber} and~\ref{t:fiberconv} imply: 
\begin{cor}\label{c:three}
Let $M$ be a closed oriented $3$-manifold with $b_1(M)=1$. 
Then $M$ is geometrically formal if and only if it fibers over $S^1$.
\end{cor}
This includes many non-symmetric manifolds.

Thurston has proved that every $3$-manifold which fibers over $S^1$
carries a unique locally homogeneous geometry. It is not clear whether
the induced metric is formal, even when the first Betti number is one.

\section{Dimension four}\label{s:four}

If $M$ is a closed oriented geometrically formal $4$-manifold, then 
by Theorem~\ref{t:elem} we have $b_1(M)\in\{ 0,1,2,4\}$.
If the first Betti number is maximal, then Theorem~\ref{t:fiber} says 
that $M$ is the 4-torus.

\subsection{First Betti number $=2$.}
In this case the Euler characteristic vanishes by Theorem~\ref{t:Euler}, 
and $b_2(M)=2$. 

Fix a formal Riemannian metric $g$. If $\alpha$ and
$\beta$ are harmonic $1$-forms generating $H^1(M)$, then they are
pointwise linearly independent. Therefore $\omega=\alpha\wedge\beta$
is a non-zero harmonic $2$-form with square zero. Thus the intersection
form of $M$ is indefinite, and we conclude $b_2^+=b_2^-=1$.
This means that the real cohomology ring of $M$ is the same
as that of the globally symmetric space $S^2\times T^2$. 

We know from Theorem~\ref{t:fiber} that $M$ fibers over $T^2$. The
above discussion shows that the fiber is nontrivial in homology.

There are many examples of such manifolds, other than $S^2\times T^2$. 
If $N$ is any $3$-manifold with $b_1(N)=1$ which fibers over the circle,
then the product $M=N\times S^1$ is a $4$-manifold with the real 
cohomology ring of $S^2\times T^2$. By Corollary~\ref{c:three} it is 
geometrically formal, because we can take a product metric which on $N$ 
is the formal metric constructed in the proof of Theorem~\ref{t:fiberconv}.


\subsection{First Betti number $=1$.}
If the first Betti number is one, the Euler characteristic vanishes
by Theorem~\ref{t:Euler}, and therefore $b_2(M)=0$. In this case
$M$ has the real cohomology algebra of the globally symmetric
space $S^3\times S^1$. 
Theorems~\ref{t:fiber} and~\ref{t:fiberconv} imply: 
\begin{cor}\label{c:four}
Let $M$ be a closed oriented $4$-manifold with $b_1(M)=1$ and $b_2(M)=0$. 
Then $M$ is geometrically formal if and only if it fibers over $S^1$.
\end{cor}
This includes many non-symmetric manifolds. The simplest example is a 
product of $S^1$ with a rational homology $3$-sphere which is not symmetric.

\subsection{First Betti number $=0$.}
From Theorem~\ref{t:elem} we know $b_2^{\pm}\leq 3$. If $b_2^+>0$, 
then there are nontrivial self-dual harmonic forms. By~\eqref{eq:first} 
they have no zeros and so define almost complex structures compatible 
with the orientation of $M$. It follows that $b_2^+$ is odd. Similarly, 
if $b_2^- > 0$, then there are almost complex structures compatible 
with the orientation of $\overline{M}$ and $b_2^-$ must be odd. 

Suppose now that $b_2^+(M)=3$. Then the self-dual harmonic forms
trivialise $\Lambda^2_+$, and each defines an almost complex structure
with trivial first Chern class (because the pointwise orthogonal 
complement of each in $\Lambda^2_+$ is trivial). Thus 
$0=c_1^2(M)=2\chi (M)+3\sigma (M)=4+5b_2^+ -b_2^- = 19-b_2^-$, 
which contradicts $b_2^-\leq 3$. Therefore $b_2^+=3$ is not possible, 
and similarly $b_2^-=3$ is not possible either.

Thus the only possible values for $b_2^{\pm}$ are $0$ and $1$, and all combinations occur for the globally symmetric spaces $S^4$, $\C P^2$,
$\overline{\C P^2}$ and $S^2\times S^2$. This finally completes
the proof of Theorem~\ref{t:main} in the $4$-dimensional case.

Any other example must have the same real cohomology ring as one of 
the above. In fact, other examples exist for each cohomological type. 
In the case of $S^4$ any rational homology $4$-sphere will do. In 
the case of $\C P^2$, there is the Mumford surface~\cite{Mu}, an 
algebraic surface of the form $\C H^2/\Gamma$ with the same rational 
cohomology as $\C P^2$. The K\"ahler form is of course harmonic, it 
generates the cohomology and its square is harmonic. Reversing the 
orientation of the Mumford surface we obtain an example with the 
cohomology ring of $\overline{\C P^2}$. Finally, in the case of 
$S^2\times S^2$, there is also a locally Hermitian symmetric 
algebraic surface $M$ of general type with the same real cohomology 
ring, due to Kuga, cf.~\cite{BPV} p.~237. In this case $M$ is of the 
form $(\HH^2\times\HH^2)/\Gamma$, and the harmonic forms are generated 
by a self-dual and an anti-self-dual harmonic $2$-form (with respect 
to the locally symmetric metric). These are K\"ahler forms for $M$ 
and $\overline{M}$ respectively, and are therefore parallel and have 
harmonic products. 

\section{Obstructions from symplectic geometry}\label{s:sympl}

In this section we discuss relations between harmonicity of products 
of harmonic forms on $4$-manifolds on the one hand, and existence
of symplectic structures on the other. This leads to further
obstructions to geometric formality. 

Let $M$ denote a closed oriented $4$-manifold with a Riemannian 
metric $g$. Suppose that $b_2^+(M)>0$. Then there is a nontrivial 
$g$-self-dual harmonic $2$-form $\omega$. If the product
$$
\omega\wedge\omega = \omega\wedge *\omega = \vert\omega\vert^2dvol_g
$$
is harmonic, then $\omega$ has constant length, and in particular has
no zeros. It is then a symplectic form on $M$ compatible with the 
orientation, and $g$ is an almost K\"ahler metric.

There are $4$-manifolds for which the elementary obstructions
of Section~\ref{s:elem} vanish, but which are not geometrically
formal because they do not admit any symplectic structure: 
\begin{ex}
Let $X$ be $\C P^2$, $S^2\times S^2$, or the Kuga or Mumford surface.
Let $N$ be a rational homology $4$-sphere whose fundamental group has 
a nontrivial finite quotient. Then $M=X\# N$ has the real cohomology 
ring of the geometrically formal manifold $X$, but is not itself
geometrically formal because it does not admit any symplectic
structure by the result of~\cite{KMT}.
\end{ex}

There is another application of the relationship between 
harmonicity of products of harmonic forms and symplectic 
structures. Namely we can show that on certain manifolds
{\it all} products of certain harmonic forms are non-harmonic.
This is obviously much stronger than geometric non-formality.

For an example, consider the smooth manifold $M$ underlying
a complex K3 surface. Then $M$ is simply connected with $b_2^+=3$
and $b_2^-=19$. The elementary considerations in Section~\ref{s:four} 
already show that $M$ is not geometrically formal. We can sharpen 
this as follows: 
\begin{prop}\label{p:K3}
Let $g$ be an arbitrary Riemannian metric on the K3 surface $M$. 
If $\alpha$ is a $g$-anti-self-dual harmonic $2$-form, then it must
have a zero. In particular the wedge product $\alpha\wedge\beta$ is 
not harmonic for any $\beta$ unless it vanishes identically. For 
example, if $\alpha$ is nontrivial then $\alpha\wedge\alpha$ is not 
harmonic.
\end{prop}
\begin{proof}
Suppose $\alpha$ is nontrivial and anti-self-dual. We have
$$
\alpha\wedge\alpha = - \alpha\wedge *\alpha = 
- \vert\alpha\vert^2 dvol_g \ ,
$$
which is harmonic if and only if the norm of $\alpha$ is constant. If it 
is constant, it must be a non-zero constant, and then the above equation
shows that $\alpha$ is a symplectic form inducing the opposite (non-complex)
orientation on $M$. In particular, $\overline{M}$ must
have non-trivial Seiberg-Witten invariants, see~\cite{T1}.
But the K3 surface contains smoothly embedded (-2)-spheres, which become 
(+2)-spheres when the orientation is reversed, showing that all the 
Seiberg-Witten invariants vanish, see~\cite{orient}.
\end{proof}
\begin{rem}
The vanishing of the Seiberg-Witten invariants for $\overline{M}$ can 
also be proved without appealing to the vanishing theorem for spheres 
of positive self-intersection. For a scalar-flat Calabi-Yau metric the 
Seiberg-Witten equations on $\overline{M}$ have no solution, though they
do have a (unique) solution on $M$. 
\end{rem}

This can be generalised quite substantially.
If $M$ has an indefinite intersection form, there are both
self-dual and anti-self-dual harmonic forms for all metrics. 
If the square of such a form is harmonic, it is a symplectic
form on $M$, respectively $\overline{M}$. But by the results 
of~\cite{orient}, manifolds which are symplectic
for both choices of orientation are quite rare.
Thus, the above proof generalises to many cases to show that
for all metrics on certain $4$-manifolds, all self-dual and/or
all anti-self-dual harmonic forms must have zeros and
non-harmonic squares. This generalises the existence of zeros
of harmonic $1$-forms on surfaces, cf.~Section~\ref{s:two}.

In the case of complex surfaces, Theorem 1 of~\cite{orient} 
implies the following:
\begin{thm}
Let $M$ be a compact complex surface of general type. Assume
one of the following conditions holds:
\begin{enumerate}
\item $K_M$ is not ample, or
\item $c_1^2(M)$ is odd, or
\item the signature $\sigma (M)$ is negative, or is zero and 
$M$ is not uniformised by the polydisk.
\end{enumerate}
Then for every Riemannian metric $g$ on $M$, all $g$-anti-self-dual
harmonic $2$-forms have zeros and non-harmonic squares.
\end{thm}
In the first case, the argument is the same as in the proof
of Proposition~\ref{p:K3}, because ampleness of $K_M$ only
fails if there are rational curves of negative self-intersection.
Note that we only need smoothly rather than holomorphically
embedded spheres, so one can replace condition 1.~by a 
weaker assumption.

\section{General existence of non-formal metrics}\label{s:generic}

Having seen that only very few manifolds are geometrically formal, 
we now want to show that even these tend to also have non-formal 
metrics. The two-dimensional case of the following result was 
already proved in Section~\ref{s:two}.

\begin{thm}\label{t:generic}
A closed oriented manifold admits a non-formal Riemannian metric 
if and only if it is not a rational homology sphere.
\end{thm}
\begin{proof}
It is clear that every metric on every rational homology sphere
is formal because there are no nontrivial harmonic forms.

Conversely, assume that $M$ is a manifold with a non-zero Betti
number $b_k(M)$, for $0<k<dim(M)$. Let $g$ be a Riemannian metric
which has positive curvature operator on an open set, say a ball 
$B\subset M$, and assume it is formal.

If $\alpha$ is a nontrivial $g$-harmonic $k$-form, then 
$$
\alpha\wedge *\alpha = \vert\alpha\vert^2 dvol_g
$$
shows that $\alpha$ has constant length. Therefore, the 
Bochner-Weitzenb\"ock formula
\begin{equation}\label{eq:BW}
\Delta (\alpha ) = \nabla^*\nabla \alpha +{\mathcal R} (\alpha)
\end{equation}
for $k$-forms allows us to compute:
$$
0=\Delta (\frac{1}{2}\vert\alpha\vert^2)= g(\nabla^*\nabla\alpha,\alpha)-
\vert\nabla\alpha\vert^2 =-g({\mathcal R} (\alpha),\alpha)-
\vert\nabla\alpha\vert^2 \ .
$$
Here the term ${\mathcal R}$ is positive on $B$, because there
the curvature operator is positive. Thus $\alpha$
vanishes identically on $B$. As $\alpha$ is harmonic, the unique
continuation principle implies that $\alpha$ vanishes on all of $M$,
contradicting the assumption that $\alpha$ is nontrivial.
\end{proof}
\begin{rem}
The above proof shows that there is an open set of non-formal
metrics in the space of all Riemannian metrics (with the $C^{\infty}$
topology, say) on any manifold which is not a rational homology sphere.
\end{rem}

\bibliographystyle{amsplain}

\bigskip

\end{document}